\documentclass[a4paper]{amsart}
\usepackage{epsfig}
\usepackage{amsmath}
\usepackage{amsthm}
\usepackage{amssymb}
\usepackage{latexsym,amsfonts}
\usepackage{amscd}
\usepackage[all]{xy}
\parskip\smallskipamount
\numberwithin{equation}{section}
\newcommand{\ck}{{{k}}} 

\newcommand{\mdeg}{\mbox{deg\,}}

\newcommand{\mHom}{\mbox{Hom}\,}

\newcommand{\Pp}{{\mathbb P}}

\newcommand{\Zz}{{\Bbb Z}}

\newcommand{\cL}{{\cal{L}}}

\newcommand{\cO}{{\cal{O}}}

\newcommand{\lra}{\longrightarrow}

\newtheorem{Prop}{Proposition}[section]
\newtheorem{Thm}[Prop]{Theorem}
\newtheorem{Lemma}[Prop]{Lemma}

\let\cal\mathcal

\def\mHom{\operatorname {Hom}}
\def\mHilb{\operatorname {Hilb}}

\def\mEnd{\operatorname {End}}

\def\deg{\operatorname {deg}}

\title{Moduli of Representations, Quiver Grassmannians, and Hilbert Schemes}
\author{Lutz Hille}

\begin{document}
\maketitle

\begin{abstract}
It is a well established fact, that any projective algebraic variety is a moduli
space of representations over some finite dimensional algebra. This algebra
can be chosen in several ways. The counterpart in algebraic geometry is
tautological: every variety is its own Hilber scheme of sheaves of length one.
This holds even scheme theoretic. We use Beilinson's equivalence to get 
similar results for finite dimensional algebras, including moduli spaces and
quiver grassmannians. Moreover, we show that several already known 
results can be traced back to the Hilbert scheme construction and Beilinson's
 equivalence.
\end{abstract}


\section{Introduction}

Assume $\ck$ is an algebraically closed field and $X$ is a projective
subscheme of $\Pp^n$ defined by some homogeneous equations 
$f_1, \ldots,f_r$ in $\ck[X_0,
  \ldots, X_n]$. We want to realize $X$ as a moduli space of quiver
representations and as a quiver grassmannian in a natural way. Moreover,
we also like to have a construction making the quiver as small as possible. 
 Let $A$ be a bounded
path algebra $\ck Q/J$, where $Q$ is a finite quiver and $J$ is an ideal
of admissible relations in the path algebra $\ck Q$. Moduli spaces for quiver
representations have been defined by King in \cite{King}, a quiver
grassmannian is just the variety of all submodules of a given module
$M$ of a fixed dimension vector. We note that we can consider moduli
spaces and also quiver grassmannians with its natural scheme
structure. Moreover, any quiver grassmannian is a moduli space
(just a moduli space of submodules of a given module), and there are 
natural morphisms between quiver grassmannians
and moduli spaces. Under certain additional conditions these morphisms are
even isomorphisms. Those isomorphisms are always hidden in our construction.
Since these morphisms
can be seen explicitely in our construction we do not need any general
result for those morphisms. This is the main reason why we use line
bundles in our construction, for arbitrary vector bundles  all constructions become much
more technical. The other advantage of using line bundles is that we can 
always use modules of dimension vector $(1,1,\ldots,1,1)$ (also called thin sincere).
\medskip

\begin{Thm}\label{Thm1}
Let $X$ be any projective scheme defined by equations $f_1,\ldots,
f_r$ in a projective $n$-space $\Pp^n$. Then there exists a quiver $Q$,
an ideal $J$ in the path algebra $\ck Q$ and a $\ck Q$--module $M$ so that
\begin{enumerate}
\item 
$X$ is isomorphic to the moduli space of all indecomposable $\ck
Q/J$--modules of dimension vector $(1,\ldots,1)$, and \\
\item 
$X$ is isomorphic to the quiver grassmannian of all submodules of
$M$ of dimension vector $(1, \ldots,1)$ for the quiver $Q$.
\end{enumerate}
\end{Thm}

Note that $Q$ can be chosen to be the Beilinson quiver, $J$ is an
ideal just defined by the $f_i$ and $M$ is the unique sincere
injective cover of one simple module, in particular $M$ is
indecomposable. For more details we refer to section \ref{sect2}. The result
on quiver grassmannians recently attracted attention in connection
with cluster algebras (see  \cite{KellerSch}) and in
connection with Auslanders theory on morphisms determined by
objects. Ringel has already pointed out that the result above has been
studied by several authors (\cite{Ringel1}), however 'can be traced
back to Beilinson' (\cite{Ringel2}). The principal aim of this note is
to show how, we can use Beilinson, and even better, how we can even improve
it. Eventually, we show that all the constructions at the end can be
traced back to a tautological construction in algebraic geometry. Any
scheme $X$ is its own Hilbert scheme of sheaves of lenght one:
$\mHilb^1(X) = X$.
\medskip

We note that
the second result, for a variety $X$, was already  stated in
\cite{H-Z2} and proven 
again with different methods in \cite{Reineke}. It can certainly be
traced back to the work in \cite{BH-Z1,BH-Z2}.  An affine version was
already proven in \cite{H-Z1}. However, the first published result in this
direction was just an example
in  \cite{HilleTilt}. We will give a common frame for all those examples, in fact 
all are variants of the Hilbert scheme construction in algebraic geometry and
a variant of Beilinsons equvivalence. Ringel already noticed that we can even
work with the Kronecker quiver, thus, two vertices are sufficient for $Q$. 
Improving this construction slightly, we can 
even realize any projective subscheme of the $n$--dimensional projective
space as a quiver grassmannian for the $(n+1)$--Kronecker
quiver. This construction is again explicit. We denote the the $m$th 
homogenous component of the ideal $I$ generated by the polynomials 
$f_i$ by $I_m$. Thus the homogeneous coordinate ring of $X$ is just 
$\oplus S^mV/I_m$ for some $(n+1)$--dimensional vector space $V$. 
We also denote by $d$ a natural number greater or equal to the 
maximal degree of the the polynomials $f_i$.
\medskip

\begin{Thm}\label{ThmKronecker}
Let $X$ be any projective subscheme of the $n$--dimensional projective space 
$\Pp^n$. There is a module $M = (S^{d-1}V/I_{d-1}, S^dV/I_d)$ over the Kronecker
algebra defined by the natural map $S^{d-1}V/I_{d-1} \otimes V \lra  S^dV/I_d$. Then
$X$ is isomorphic to the quiver grassmannian of submodules of $M$ of dimension 
vector $(1,1)$.
\end{Thm}

The principal part of the note consists of a five step construction
that we will use to get a realisation of $X$ as such a moduli
space. In addition we also add some modifications of these steps 
allowing to simplify the quiver or the relations.
 We explain these five steps briefly. First note, that any scheme
$X$ is its own Hilbert scheme of sheaves of length one. So any
projective scheme is a moduli space of sheaves (in a rather trivial
way). In a second step we use Beilinson's equivalence to construct for
$X$ an algebra $A = \ck Q /J$. Roughly, we can take any tilting bundle
$T$ on $\Pp$, extend it by any other vector bundle $T'$ to $R = T \oplus T'$
and apply $\mHom(R,-)$ to the universal family of the Hilbert scheme
$X$. In the particular case when $T$ is the direct sum of the line
bundles $\cO(i)$, for $i = 0, \ldots, n$, we can just extend it by the
line bundles $\cO(i)$ for $i=n+1, \ldots,d$. In this way, we get a
family of modules over the Beilinson algebra for $T$ and a family of
modules over the 'enlarged' Beilinson algebra for $R$. If $X$ is given
by polynomials as above, all modules of the family also satisfy the
equations $f_i$, however now in the Beilinson algebra. Thus we define
$J$ to be the ideal generated by the $f_i$ in the enlarged Beilinson
algebra. Note that we have several choices for such realization,
depending on where the relation starts. However, one can check
directly, that family of all modules of dimension vector
$(1,\ldots,1)$ over $A$, the enlarged Beilinson algebra with relations $J$,
coincides with $X$, independent of the this realization. Consequently, 
the moduli space of all modules of
dimension vector $(1, \ldots,1)$ over $A$ is $X$ as a scheme. In a
final step, we realize $X$ as a quiver grassmannian by using an injectice
hull in $A$.
\medskip

We already mentioned that this construction is more general in the way
that we can take any tilting bundle $T$ and any vector bundle $T'$,
however the direct computation seems to be more sophisticated. So we
use line bundles just for simplicity. Even stronger, in general we do
not need $T$ to be a tilting bundle. For example, the construction
also would work if we only take $\cO \oplus \cO(d)$ where $d$ is at
least the maximal degree of the $f_i$. Then we get a realization
similar to Reinekes construction, that is in fact a variant of the 
realization for the Kronecker quiver.
\medskip

The paper is organized as follows. In the second section we present
the five steps of our construction together with a final note on
framed moduli spaces. In the third section we reprove the already
known results using the construction in section \ref{sect2}. We
conclude in the last section with some open problems and a proof of 
the second theorem.

\medskip



\section{Hilbert schemes, Beilinson's equivalence and Serre's
  construction}\label{sect2} 

We construct, using some elementary results from algebraic geometry,
for any algebraic variety a moduli space of quiver representations, a
quiver grassmannian and also further examples in five steps. We also
note, that this even holds for any scheme that is quasi-projective. So
we obtain the authors example from 1996 \cite{HilleTilt}, Huisgen-Zimmermann's
examples in her work on uniserial modules (see for example
\cite{H-Z1}, and her work with Bongartz \cite{BH-Z1, BH-Z2} on
Grassmannians, we apologize 
for being not complete), a variant of Reineke's result
for quiver-grassmannians from 2012 \cite{Reineke} and last not least
Michel Van den Bergh's example, that appeared in a blog of Lieven le
Bruyn. 
In fact, all the
results at the end of this note can be  
proven using the following constructions.

\subsection{Hilbert schemes} (\cite{HilbertSchemes})

 We take  an algebraic variety $X$ and consider sheaves of length one
 on $X$. There is a bijection between those sheaves and
 points of $X$. In a more sophisticated way we can say $X =
 \mHilb^1(X)$ the Hilbert scheme of length one sheaves on $X$. Or, we
 can consider any line bundle $\cL$ on $X$ as a fine moduli space of
 skyscraper sheaves by taking the push forward to the diagonal in $
 X \times X$. Each fiber of some point $x$ for the first projection is
 just the skyscraper sheaf in $x$.
 
 \subsection{Beilinson's tilting bundle} (\cite{Beilinson})
 
 In the second step we transform the construction above to the
 representation theoretic side using tilting. 
 To keep the construction easy, we consider the Beilinson tilting
 bundle $T = \cO \oplus \cO(1) \oplus ... \oplus \cO(n)$ on the  
 $n$--dimensional projective space $\Pp^n$. We denote by $A$ the
 Beilinson algebra, it is the opposite of $\mEnd(T)$. Take $\ck_x$ 
  to be a skyscraper sheaf and let us compute $\mHom(T,\ck_x)$. Note
  that  $\dim \mHom(\cO(i),\ck_x) = \ck $, thus we get thin sincere  
  representations of $A$, that is each simple occurs with multiplicity
  one. Now we can see by direct computations  that the moduli  
  space of  thin sincere representations of $A$ is just the projective
  space $\Pp^n$. 
 
 This example can be easily generalized to any tilting bundle $T$,
 however, using sheaves (that are not vector bundles) we do not get a
 flat family (the dimension above jumps at certain points). 
 
\subsection{Relations of length at most $n$} (\cite{Serre})

Now we consider $Y$, any subvariety defined by equations $f_1, \ldots,
f_r$ in a projective $n$--space $\Pp^n$. Assume first $\deg f_i \leq
n$ for all $i$. Note that the quiver of the Beilinson algebra above
has $n+1$ arrows from $i$ to $i+1$, we denote by $x_0^i,\ldots,
x_n^i$, and paths the monomials
$x_{i(1)}^ax_{i(2)}^{a+1}...x_{i(d-1)}^{d-1}x_{i(d)}^d$ of degree
$d-a+1$ at most
$n$. We define an 
algebra $B$ as the quotient of $A$ by the ideal $J= (\overline f_1,
\overline f_2, \ldots, \overline f_r)$, where $\overline f$ is any
linear combination of path representing $f$ in $A$. Note that any
representative works, since the arrows in $A$ commute, whenever this
makes sense: $x_i^ax_j^{a+1} = x_j^ax_i^{a+1}$.  

This solves the problem we address to the next step already,
since $\Pp^n \subset \Pp^m$ for any $n < m$. Moreover, we could also
use the fact, that any projective algebraic variety is already defined
by quadratic relations. However, using Serre's theorem we can handle
also relations of any degree in $\Pp^n$. 
\medskip

\begin{figure}[ht]
	\centering
\includegraphics[width=10cm]{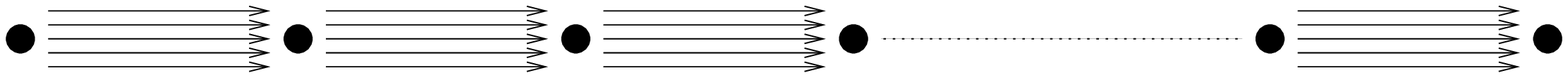}
\caption[Bildunterschrift]{Beilinson quiver}
	\label{Beilinson quiver}
\end{figure}
\medskip

\subsection{Relations of arbitrary length} ~

To obtain $Y$, as in the previous step, where $\mdeg f_i \leq d$ for
any $d > n$, we consider the sequence of line bundles $\cO, \ldots,
\cO(d)$. The 
direct sum of these line
bundles is no longer a tilting bundle, however the same computation as
above shows that the moduli space of all thin sincere representations
of $A$ is still a projective $n$--space. Now the representatives
$\overline f_i$ of
the polynomials $f_i$ live in $A = \mEnd(\oplus_{i=0}^d \cO(i))$ and the
moduli space of thin sincere representations of $B = A/(\overline f_1,
\ldots, \overline f_r)$ is $Y$ (even scheme theoretic). The reader
familiar with Serre's construction will notice that this step is just
inspired by this construction (\cite{Serre}). 

\subsection{Quiver Grassmannians} ~

In a final step, we use the embedding of the thin sincere
representations of $B$ in its minimal injective hull $M$. Note that
$I$ is the indecomposable injective module with the unique simple 
socle, that is the socle
of any thin sincere representation of $B$. In case $B = A$ (all $f_i$
are zero) we obtain the projective space as quiver grassmannian of
thin sincere subrepresentations of $I$. In a similar way, also $Y$
coincides  with the quiver grassmannian of thin sincere 
subrepresentations of the large indecomposable injective $B$--module
$M$. This final step goes back to Schofield \cite{Schofield1} and was
mentioned later 
also by Van den Bergh and Ringel (\cite{Ringel2}).

\subsection{Framed moduli spaces}~

We note that quiver grassmannians can be obtained also directly from a
corresponding construction in algebraic geometry. Any
skyscraper sheaf is the quotient of a line bundle $\cL \lra k_x$. If
we consider the moduli space of all those quotients of $\cL$ with fixed
Hilbert series of the quotient sheaf, we get the (framed) Hilbert scheme,
that coincides with the original one. If we apply Beilinson's tilting  again, we
get for $\cL = \cO(d)$ a projective $A$--module. Thus $\Pp^n$ is the
quotient grassmannian for the large (that is sincere) indecomposable
projective $A$--module. The same construction works with $Y$ instead
of $\Pp^n$.

\subsection{Reduction of the quiver} ~

Note that any vector bundle $T = \oplus_{i \in L} \cO(i)$ for any $L$ with at least two
elements on  $\Pp^n$ defines a morphism from 
$\Pp^n$ to the moduli space of modules of dimension vector $(1,1,\ldots,1,1)$ 
over $A = \mEnd(T)$ (or even over the path algebra of $A$) and also to the 
corresponding quiver grassmannian for the injective hull $M$ of an indecomposable
of dimension vector $(1,1,\ldots,1,1)$. Thus, for the projective $n$--space even two
line bundles are sufficient. In the last section we modify this construction slightly and
consider $L$ consisting of three, repectively even two, elements so that we still get 
an isomorphism for any subscheme $X$ in $\Pp^n$. 

\subsection{Proof. }
We prove the  Theorem \ref{Thm1} using the five steps above. In explicit terms
the module $M$ is just defined by vector spaces $M_m = S^mV/I_m$ the 
$m$th graded part of the homogeneous coordinate ring. This becomes a 
module over the Beilinson quiver using the natural map 
$V \otimes S^mV \lra S^{m+1}V$ as a multiplication map as follows. Take a 
basis $v_0, \ldots, v_n$ of $V$ and define the linear map of the $i$th arrow
just by tensoring with $v_i$: $S^mV \lra S^{m+1}V$. The commutative relations 
force that the moduli space (or the corresponding quiver grassmannian in $M$)
of thin sincere modules is just a subscheme in $\Pp^n$ defined by some of the 
polynomials $f_i$. If we consider sufficiently many degrees $m$, then any $f_i$ is
realized in some homogeneous part $I_m$ of the ideal $I$. For example the two 
degrees $m = 0, d$ are sufficient to see any $f_i$. Thus, if we take the Beilinson 
quiver with vertices $0, 1, 2,  \ldots, d-1, d$ we can certainly realize the variety $X$.
\medskip

Conversely, we may ask how many degrees we need to realize $X$ in the module $M$.
A similar consideration as above shows that the three degrees $m = 0,e,d$ are sufficient, 
provided $d$ is at least the maximum of all degrees of the polynomials $f_i$, and $e$ 
can be any natural number with $0 < e < d$. In particular, we can take $e=1$ or $e= d-1$.
This leads to the proof of the second theorem proven in the last section. The lemma 
below reduces than even to the Kronecker quiver.

\section{Overview on results and some consequences}

\subsection{Some variatians} 

Now we can use the construction above to get many variations, we can
not list all, however we should collect some. First we 
construct affine examples. One way is to take open subvarieties,
however we would like to characterize open subsets module theoretic.
Huisgen-Zimmermann started to consider uniserial modules in
\cite{H-Z1}. To obtain 
affine varieties as moduli spaces of uniserial modules, we consider a 
variant of the Beilinson quiver, we replace the first arrow $x_0$ just
by a path $y_0z_0$ of length two. Then a thin sincere module is
uniserial precisely when its map $y_0z_0$ is not zero. Thus it is the
open subvariety (subscheme) defined by $x_0 = 1$, that is an affine
chart.  

\begin{figure}[ht]
	\centering
\includegraphics[width=10cm]{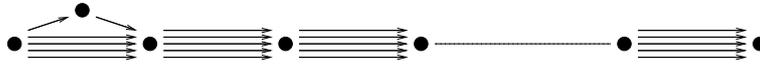}
\caption[Bildunterschrift]{modified Beilinson quiver}
	\label{modified Beilinson quiver}
\end{figure}

\subsection{Consequences} ~

In our opinion there are two kind of consequences. First one might
think that we can now obtain results in algebraic geometry using
representation theory. This seems to be impossible, as far we
consider any algebraic variety.  However, restricting to some
subclasses this might be fruitful, we mention some open problems 
at the end of this note. Moreover, for our construction,
using the Beilinson algebra, the relations are directly given by the
defining polynomials. Thus we do not get any deeper insight by
considering an algebraic variety as a moduli space of quiver
representations. 

The second consequence concerns the realisation of a variety as a particular
moduli space, that is more restrictive. This is often very  useful and
is already used quite often. The main open problem here seems to be
to construct all moduli spaces of quiver representations for a particular
quiver. In general, for all dimension vectors, this is even open for the
$3$--arrow Kronecker quiver.

\subsection{Results}

We use the construction in the previous section to prove some of the
already known results just by applying the five steps. We start with
any projective algebraic variety and proceed with affine ones. As we
already explained, we consider $X$ as the scheme of length one sheaves
on itself and apply the Beilinson tilting bundle.

\begin{Thm} \cite{HilleTilt}
Any projective algebraic scheme of finite type is a fine moduli space
of modules over some finite dimensional algebra (a bounded path
algebra). Moreover, we can 
obtain it already for the thin sincere representations, that is the
Jordan-H\"older series contains each simple module just once in its
composition series up to isomorphism. 
\end{Thm} 

Taking open parts we recover the result of Huisgen-Zimmermann, that was
obtained using uniserial modules
(Theorem G in \cite{H-Z1}). Note, the result was stated in \cite{H-Z1} in a different 
language, the notion of a moduli space was adapted by her only later.


\begin{Thm} \cite{H-Z1}
Any affine algebraic variety 
is a fine moduli space of uniserial
modules over some finite dimensional algebra (a bounded path
algebra). 
\end{Thm}

Then Grassmannians also appeared in Huisgen-Zimmermann's work, however
the idea was 
already introduced by Schofield \cite{Schofield1} and then intensively used by
Nakajima \cite{Nakajima}. However, a similar result could be read of
from the work of Bongartz and 
Huisgen-Zimmermann and was later explicitely stated in \cite{H-Z2}. Here
again we can use thin sincere submodules of a module $M$ or just uniserial 
modules.

\begin{Thm} \cite{H-Z2} 
Any projective algebraic variety is a quiver grassmannian.
\end{Thm} 


\subsection{Kronecker quiver}
Using the Beilinson construction with a rather small vector bundle we
can reduce the quiver even to the Kronecker quiver. This is almost
the same construction as in Reinekes work and based on the following
geometric construction. Consider $\Pp^n = \Pp(V)$ embedded into 
$\Pp(S^mV)$ with the $m$--uple embedding. Assume $X$ is a subscheme
in $\Pp(V)$ and consider its image in $\Pp(S^mV)$. If $m$ is larger
than the maximal degree of the polynomials $f_i$ defining $X$, the 
equations of $X$ in $\Pp(S^mV)$ are just linear and the defining equations of the
embedding $\Pp(V) \lra \Pp(S^mV)$ (that are quadratic). Just modilfying 
the Beilinson construction we can use the bundle $\cO \oplus 
\cO(e) \oplus \cO(d)$. This reduces the construction to a quiver with
three vertices. For $e = d/2$ and $d$ sufficiently large, this corresponds
to realizing $X$ using quadratic equations. The corresponding module $M$
considered as a representation of a three vertex quiver $(M_1, M_2, M_3)$
has a simple socle $M_3 = k = S^0V$ with $M_2 = S^eV$ and $M_3 = S^dV$. 
Now we use Ringels idea to reduce to the Kronecker quiver $S^{d-e}V$.

\begin{Lemma}
With notation above and any $d > e > 0$ we have an isomorphism of quiver 
grassmannians as follows. 
The quiver grassmannian of submodules of $M = (M_1, M_2, M_3)$ of 
dimension vector $(1,1,1)$ 
is isomorphic to the quiver grassmannian of submodules of $(M_1,M_2)$
of dimension vector $(1,1)$.
\end{Lemma}

{\sc Proof. }
Note that the restriction of $M$ to $(M_1, M_2)$ defines a morphism of
quiver grassmannians. Since $M_3$ is just one--dimensional, any submodule
$(M_1, M_2)$ over the Kronecker algebra of dimension vector $(1,1)$ 
extends uniquely to a submodule of $M$ of dimension vector $(1,1,1)$. Thus,
this morphism is a bijection. In the particular case of $X$ being the projective
space, this morphism is an isomorphism. Going back to $M$ we just restrict
this isomorphism to the subscheme defined by the polynomials $f_i$, 
consequently, both quiver grassmannians are also isomorphic.
\hfill $\Box$
\medskip

Taking $d$ at least the degree of the defining equations $f_i$ and $d = e+1$
we realize $X$ as a quiver grassmannian over the Kronecker algebra with $(n+1)$
vertices. This proves Theorem \ref{ThmKronecker}. Note that Reineke realized
the linear subspace using an additional arrow, however, this is not necessary.
\medskip

Obviously, we can not reduce to just one vertex, thus two vertices is the
minimum we can achieve. However, it is not clear whether we can
still reduce the number of arrows. Such a reduction would be more
complicated and certainly independent of Beilinsons result.

\subsection{Affine versus projective} 
At the end we discuss the problem how to obtain projective examples
from affine ones and vice versa. As we have mentioned above, one can
take a projective variety, that is a moduli space, and obtain an
affine cover as moduli spaces of uniserials by modifying the Beilinson
quiver slightly.

The converse, to obtain complete examples by glueing, is an open
problem. In particular, let $X$ be a complete variety that is not
projective (see Hartshorne for an example \cite{Hartshorne}, Ex 3.4.1
in appendix B)
then to our knowledge there is no way so far, to get $X$
as a moduli space of representations. Moreover, it is clear that $X$
can not be a quiver grassmannian, since the latter one is projective by
definition. One might think that also moduli spaces are always projective,
however, we should mention that moduli spaces as constructed in King's
paper \cite{King} are, but there might be other constructions as well.

\subsection{Further open problems}
Since already Kronecker quivers are very complicated with respect to
the geometry of quiver grassmannians it would be natural to restrict
to particular classes of modules or quivers. As far we know, the problem 
to describe all quiver grassmannians is open for Dynkin quivers and also 
tame quivers. It also would be desirable to understand quiver grassmannians
for the $3$--arrow Kronecker quiver. Moreover, inspired by cluster algebras,
the main open problem seems to be to understand quiver grassmannians
for exceptional modules over path algebras.
\medskip

If we use the explicit construction of the module $M$ with $M_m = S^mV/I_m$ one
can see, that everything is even defined over any base field. For polynomials 
over the integers everything is defined even over $\Zz$. Thus the construction
also works in the same fashion over an commutative ring. 




\begin{thebibliography}{References}
\setlength{\parskip}{0em}


\bibitem{AltmannHille}{Altmann, Klaus; Hille, Lutz: {\em Strong
    exceptional sequences provided by quivers.}
  Algebr.~Represent.~Theory 2 (1999), no.~1, 1–17.} 

\bibitem{Beilinson}{Beĭlinson, A. A.: {\em Coherent sheaves on $P^n$
    and problems in linear algebra. } (Russian) Funktsional. Anal. i
  Prilozhen. 12 (1978), no. 3, 68–69. } 

\bibitem{BGG}{Bernstein, I. N.; Gelʹfand, I. M.; Gelʹfand, S. I.: {\em
    Algebraic vector bundles on $P^n$ and problems of linear
    algebra. } (Russian) Funktsional.~Anal.~i Prilozhen.~12 (1978),
  no.3, 66–67.}    

\bibitem{BH-Z1}{ Bongartz, Klaus; Huisgen-Zimmermann, Birge: {\em  The
  geometry of uniserial representations of algebras. II. Alternate
  viewpoints and uniqueness.} J.~Pure Appl.~Algebra 157 (2001), no.~1,
  23–32.}  

\bibitem{BH-Z2}{ Bongartz, Klaus; Huisgen-Zimmermann, Birge: {\em
    Varieties of uniserial representations. IV. Kinship to geometric
    quotients.} Trans.~Amer.~Math.~Soc.~353 (2001), no.~5, 2091–2113. }  
 

\bibitem{Hartshorne}{Hartshorne, Robin: {\em Algebraic Geometry.}
  Graduate Texts in Mathematics, No. 52. Springer-Verlag, New
  York-Heidelberg, 1977.} 

\bibitem{HilleTilt}{Hille, Lutz: {\em Tilting line bundles and moduli
    of thin sincere representations of quivers.}  Representation
  theory of groups, algebras, and orders (Constanţa,
  1995). An.~Ştiinţ.~Univ.~Ovidius Constanţa Ser.~Mat.~4 (1996),
  no.~2, 76–82.}


\bibitem{H-Z1}{Huisgen-Zimmermann, Birge: {\em The geometry of
    uniserial representations of finite-dimensional algebra. I.}
  J.~Pure Appl.~Algebra 127 (1998), no.~1, 39–72.} 

\bibitem{H-Z2}{Huisgen-Zimmermann, Birge: {\em Classifying
    representations by way of Grassmannians.}
  Trans. Amer. Math. Soc. 359 (2007), no. 6, 2687–2719.} 

\bibitem{KellerSch}{Keller,  Bernhard; Scherotzke,  Sarah: 
  {\em  Desingularizations of quiver Grassmannians via graded quiver
    varieties. }
 arXiv:1305.7502. }
   

\bibitem{King}{King, Alastair: {\em  Moduli of representations of
    finite-dimensional algebras. } Quart.~J.~Math.~Oxford Ser.~(2) 45
  (1994), no.~180, 515–530.} 

\bibitem{HilbertSchemes}{Nakajima, Hiraku: {\em  Lectures on Hilbert
    schemes of points on surfaces.}  University Lecture Series,
  18. American Mathematical Society, Providence, RI, 1999.}

\bibitem{Nakajima}{Nakajima, Hiraku: {\em Quiver varieties and
    Kac-Moody algebras.}  Duke Math.~J.~91 (1998), no.~3, 515–560.}

\bibitem{Reineke}{Reineke, Markus: {\em 
    Every projective variety is a quiver Grassmannian. }
    arXiv:1204.5730.}

\bibitem{ReinekeChow}{Reineke, Markus: {\em The Harder-Narasimhan
    system in quantum groups and cohomology of quiver moduli. }
  Invent.~Math.~152 (2003), no.~2, 349–368.}  

\bibitem{Ringel1}{Ringel, Claus Michael: {\em The Auslander bijections:
    How morphisms are determined by modules. } arXiv:1301.1251.} 

\bibitem{Ringel2}{Ringel, Claus Michael:{\em Quiver Grassmannians and
    Auslander varieties for wild algebras. } arXiv:1305.4003. }

\bibitem{Schofield1}{Schofield, Aidan: {\em  General representations
    of quivers.} Proc. London Math. Soc. (3) 65 (1992), no. 1, 46–64.}

\bibitem{Schofield2}{Schofield, Aidan: {\em Birational classification
    of moduli spaces of representations of quivers.}
  Indag.~Math.~(N.S.) 12 (2001), no.~3, 407–432.}  

\bibitem{Serre}{Serre, Jean-Pierre: {\em Faisceaux algébriques
    cohérents.} (French) Ann.~of Math.~(2) 61, (1955). 197–278. }  


\end{thebibliography}
\end{document}